\documentclass[a4paper,12pt]{article}
\usepackage{amssymb}\usepackage{amsmath}\usepackage{epsfig}

\def\halfbarlength{12 pt}
\def\node{{\kern -1pt\circ\kern -1pt}}
\def\bulnode{{\kern -1pt\bullet\kern -1pt}}
\def\dbulnode#1{{\kern-1pt\displaystyle\mathop\bullet_{\hbox to 0pt{\hss
$\scriptstyle#1$\hss}}\kern-1pt}}

\def\nbar#1{{\vrule width\halfbarlength height2.75pt depth-2.25pt}_{\hbox to -1pt{\hss
$\scriptstyle _{#1}$\hss}}{\vrule width\halfbarlength height2.75pt depth-2.25pt}}
\def\vertbar#1#2{\rlap{\kern 0.7pt\vrule width0.5pt height1pt depth16pt}
                 \rlap{\lower16.7pt\rlap{$#2$}}#1}
\def\dvertbar#1#2#3{\rlap{\kern 0.7pt\vrule width0.5pt height16pt

depth1pt}\rlap{\kern 0.7pt\vrule width0.5pt height1pt depth16pt}
                 \rlap{\lower -16.7pt\rlap{$#2$}}\rlap{\lower16.7pt\rlap{$#3$}}#1}

\newbox\proofbox
\def\proof{\futurelet\next\lookforbracket}
\def\lookforbracket{\ifx\next[\let\go\usespecialterm
\else\let\go\relax
\ifvmode\vskip-\lastskip\fi\global\setbox\proofbox=\vbox
\bgroup%
\noindent{\\\\\sc Proof~:}%
\enskip\relax\fi\ignorespaces\go
}
\def\usespecialterm[#1]{\ifvmode\vskip-\lastskip\fi
\global\setbox\proofbox=\vbox\bgroup%
\noindent\hskip\parindent%
{\\\noindent\sc Proof #1:}\textnormal\ \ \relax
\ignorespaces}
\def\endproof{\hbadness10000\parfillskip=0pt\egroup
\unvbox\proofbox
\setbox0=\lastbox
\ifdim\ht0>1pt
\vskip-2pt
\noindent
\hbox to\textwidth{\vbox{
\parindent=0pt
\unhbox0 $\square$ \hss}\hss}
\relax
\fi
}

\makeatletter
\renewcommand{\@biblabel}[1]{\quad #1.}
\makeatother
\newtheorem{Lem}{Lemma}[section]

\newtheorem{Def}[Lem]{Definition}
\newtheorem{Prop}[Lem]{Proposition}

\newtheorem{The}[Lem]{Theorem}
\newtheorem{proposition}{Proposition}

\begin{document}
\title{Parabolic subgroups of Garside groups}
\author{EDDY GODELLE}
\maketitle
\begin{abstract}
\noindent A Garside monoid is a cancellative monoid with a finite lattice generating set; a Garside group is the group of fractions of a Garside monoid. The family of Garside groups contains the Artin-Tits groups of spherical type.  We generalise the well-known notion of a parabolic subgroup of an Artin-Tits group into that of a parabolic subgroup of a Garside group. We also define the more general notion of a Garside subgroup of a Garside group, which is related to the notion of LCM-homomorphisms between Artin-Tits groups. We prove that most of the properties of parabolic subgroups extend to this subgroups. keywords : Garside groups -- lattices -- parabolic subgroups.\end{abstract}
\section*{Introduction.}
The braid group on $n$ strings and more generally the Artin-Tits groups of spherical type (see Appendix A for a definition) have been the focus of many articles and are pretty well understood. It has been shown by Dehornoy and others in recent articles (see for instance \cite{DeP,Deh7,Pic2,BDM}) that many properties of Artin-Tits groups still hold for a wider class of groups, namely the Garside groups. Recall that in a monoid $G^+$ we say that $g$ left-divides ({\it resp.} right-divides) $h$ if $h = gk$ ({\it resp.} $h = kg$) for some $k$ in $G^+$. For $h$ in $G^+$, we denote by $L(h)$ and $R(h)$ the set of elements which left-divide  $h$ and right-divide $h$ respectively. We say that $G^+$ is cancellative if $ghk = gh'k$ with $g,h,h',k$ in $G^+$ implies $h = h'$ and finally we say that $G^+$ is Noetherian if for every  element $g$ in $G^+$ the set $\{k\in\mathbb{N}\mid g = g_1\dots g_k; g_i\in G^+; g_i\neq 1 \}$ is finite.
A Garside monoid is a cancellative Noetherian monoid which is a left-lattice and a right-lattice and which has an element $\Delta$ such that $R(\Delta)$ is equal to $L(\Delta)$ and is a finite generating set of $G^+$. A Garside group $G$ is the group of fractions of a Garside monoid.\\ 
One of the main properties of the Artin-Tits groups is the existence of natural subgroups, the so-called {\it standard parabolic subgroups}; see Section \ref{sec.cas.class} for a definition. Therefore it is natural to address the following question~: what are the standard parabolic subgroups of a Garside group, if they exist~?\\

The technical notion of an LCM-homomorphism (see Appendix \ref{lab.appendLCM} for a definition) was introduced in \cite{Cri} and it proves to be crucial in the solution of two long-standing conjectures on Artin-Tits groups, namely the so-called Tits conjecture and the embedding conjecture of an Artin-Tits monoid in its associated Artin-Tits group (\cite{CrP} and \cite{Par4}). A special case of LCM-homomorphism is the canonical embedding homomorphism of a standard parabolic subgroup in the group. Consider an LCM-homomorphism between Artin-Tits groups of spherical type; then this homomorphism is injective and respects the normal forms of the elements of the associated monoids and the normal forms of the elements of the groups. Hence the family of all subgroups of an Artin-Tits group of spherical type which are the image of an LCM-homomorphism is a very natural family of subgroups of the group. Thus it is natural to ask for a family of subgroups of a Garside group, which will called {\it Garside subgroup},  which extends our notion of standard parabolic subgroup and such that the image of an LCM-homomorphism between Artin-Tits groups of spherical type becomes a Garside subgroup.\\
The aim of this paper is to answer positively to the above questions.\\ 

We refer to the next sections for precise definitions of a Garside subgroup, of  parabolic subgroups and the notion of {\it sublattice, atoms, balanced element, Garside element, normal form...}. Roughly speaking, a Garside submonoid of a Garside monoid is a submonoid which is a sublattice of the monoid and which is closed by taking normal forms; that is the terms of the normal form of an element of the submonoid are in the submonoid. A Garside subgroup of a Garside group is a subgroup generated by a Garside submonoid. The two following propositions show that the nice and useful properties of parabolic subgroups of Artin-Tits groups are still true in the more general context of Garside groups; hence make our definition natural.
\begin{proposition}\label{labpropintro1} Let $G^+$ be a Garside monoid; then, \\ \noindent\begin{tabular}{cl} (a)&every Garside submonoid of $G^+$ has a Garside monoid structure;\\(b)&every element of a Garside submonoid has the same normal form\\& considered as an element of $G^+$ or considered as an element of\\& the submonoid;\\(c)&the intersection of two Garside submonoids is a Garside  submonoid.\end{tabular}\end{proposition}
\begin{proposition}\label{labpropintro2} Let $G$ be a Garside group associated to a Garside monoid~$G^+$; then,\\\noindent\begin{tabular}{c l} (d)&every Garside subgroup of $G$ is a  Garside group associated with\\&some Garside submonoid of $G^+$;\\(e)&if $H$ is a Garside subgroup of $G$ and $x$ belongs to $H$, then  the\\& fractional decomposition  of $x$ in $H$  coincides with its fractional\\& decomposition in $G$.\\ (f)&the intersection of two Garside subgroups is a Garside subgroup.\end{tabular}\end{proposition}
We say that a Garside subgroup is a {\it standard parabolic subgroup (submonoid)} when it is a subgroup (submonoid) generated by the support (in the monoid) of a balanced minimal  element of  the group and it verifies some property. In particular, its atoms are atoms of the group (monoid).  In the case of an Artin-Tits group of spherical type, with its classical presentation, {\it our} standard parabolic subgroups are the same as the classical ones; furthermore, \begin{proposition} Let $G$ be a Garside group; then, \\\noindent\begin{tabular}{c l}(g)& every standard parabolic subgroup of $G$ is a Garside subgroup of $G$.\\(h)& the intersection of two standard parabolic subgroups of $G$ is a \\&standard parabolic subgroup of $G$.\end{tabular}\end{proposition} 
We show by examples that a subgroup of a Garside group generated by a set of atoms of the group is not necessarily a standard parabolic subgroup or even a Garside subgroup. Moreover, such a subgroup can be a Garside subgroup but not be a standard parabolic subgroup. We call atomic Garside subgroups of a Garside group the Garside subgroups which atoms are atoms of the group. This family of subgroups is between the family of standard parabolic subgroups and the family of Garside subgroups. Unfortunately it is not closed by intersection as we will see in Section 1.\\  

The paper is organised as follows. In Section 1 we recall notation  and  basic results on Garside group. In Section 2 we want to explain why your definition seems to be the good one. So, firstly we recall the classical case of the Artin-Tits groups. Secondly, we study two examples in order to lead the reader to the good definitions of a Garside subgroups and of a parabolic subgroup in the context of the Garside groups. Finally in Section 3, we define the notions of a Garside subgroup and of a parabolic subgroup and then we  prove Proposition 1,2 and 3. In the appendix, we recall well-known facts, and some notation, on Artin-Tits groups which can help to read the present paper.\\

In a forthcoming paper \cite{god_prep}, we will study the normaliser of a Garside subgroup. We will introduce and investigate the category of ribbons; it will lead us to new results on Artin-Tits groups.

\section{Background and Notation for Garside group}

If $G^+$ is a monoid and $g,h$ are in $G^+$, we write $g\prec h$ and $h\succ g$ if $g$ left-divides $h$ and $g$ right-divides $h$ respectively. We say that $h$ in $G^+$ is an atom of $G^+$ if $h = gk$ with $g,k$ in $G^+$ implies $g = 1$ or $k = 1$. Recall that for $h$ in $G^+$, we set $L(h) = \{g\in G^+\mid g\prec h\}$ and $R(h) = \{g\in G^+\mid h\succ g\}$. For $g$ and $h$ in $G^+$ we denote by $g\land_L h$ and $g\lor_L h$  ({\it resp.} $g\land_R h$ and $g\lor_R h$) the g.c.d. and the l.c.m. of $g$ and $h$ for $\prec$ ({\it resp.} for $\succ$) in $G^+$. When useful to prevent confusion, we shall include $G^+$ as a subscript to specify the monoid. 

We refer to the first chapter of \cite{KoM} for the notion of lattice and to \cite{Deh7} for the general theory of Garside groups. Recall that a sublattice of a lattice is a non empty subset which is closed by $\land$ and by $\lor$. Since the first definition of Garside monoids introduced by Dehornoy and Paris in \cite{DeP}, several other definitions appeared (\cite{Bes,BDM,Deh7}). The final and generally accepted definition seems to be the following~:
\begin{Def}\label{defgarsgp}
Let $G^+$ be a monoid and $X$ a subset of $G^+$.\\
i) We say that $\Delta\in G^+$ is balanced if $R(\Delta) = L(\Delta)$. In that case, $R(\Delta)$ is denoted by $D(\Delta)$. \\ 
ii) We say that the pair $(G^+,X)$ is a   positive Garside  system if \\
\indent (a) $G^+$ is Noetherian
 and  cancellative,\\
\indent (b) $G^+$ is a lattice for right-divisibility and a lattice for left-divisibility.\\\indent (c) $X$ is a generating set for  $G^+$ and there exists $\Delta\in G^+$ balanced such that $X = D(\Delta)$.\\
If $(G^+,D(\Delta))$ is a   positive Garside  system, then $G^+$ is called a quasi-Garside monoid and $\Delta$ its Garside element. The elements of $D(\Delta)$ are called the minimals of the   positive Garside  system. When $D(\Delta)$ is finite then $G^+$ is called a Garside monoid.
\end{Def}
The quasi-Garside monoids have been recently introduced in \cite{Dig} in order to give a good divisibility structure to the affine type braid group.\\ 
If $(G^+,D(\Delta))$ is a   positive Garside  system then the set $S$ of its atoms is a subset of $D(\Delta)$.\\ Note that under the assertions (a) and (c), the assertion (b) is equivalent to \begin{center} (b') : $D(\Delta)$ is a lattice for right-divisibility and a lattice for left-divisibility,\end{center} see \cite{BDM} or \cite{Deh7} for instance.  Remark that the $\prec$ and $\succ$ are antisymmetric because the monoid is cancellative, and by the Noetherian hypothesis,  $abc = b$ implies $a = c = 1$. Also if $D(\Delta) = D(\Delta')$ then $\Delta = \Delta'$.\\
If $G^+$ is a monoid with $S$ for set of atoms and $\delta$ is a balanced element of $G^+$, then the set $D(\delta) \cap S$ is called the support of $\delta$ and is denoted  by $supp(\delta)$. 
\begin{Def}Let $G$ be a group, let $X$ be subset of $G$; denote by $G^+$ the submonoid of $G$ generated by $X$. We say that $(G,X)$ is a Garside system if $(G^+,X)$ is a   positive Garside  system and if $G$ is the group of fractions of $G^+$. In that case, we say that $G$ is a quasi-Garside group; We say that $G$ is a Garside group when $G^+$ is a Garside monoid.\end{Def}
If $(G^+,D(\Delta))$ is a   positive Garside  system  then it verifies the Ore relations. So, $G^+$ injects into its group of fractions $G$ and $(G,D(\Delta))$ is a Garside system. If $(G,D(\Delta))$ is a Garside system, we denote in the following by $G^+$ the submonoid of $G$ generated by $D(\Delta)$. In the same way that an Artin-Tits group and its Artin-Tits system are commonly confused, we will most of the time make the confusion between a Garside group ({\it resp.} monoid) and its ({\it resp.} positive) Garside system. We have the following properties:
\begin{Lem}[\cite{Deh7},\cite{Dig}]\label{decompG} Let $(G,D(\Delta))$ be a Garside system; then\\i) $D(\Delta)\Delta = \Delta D(\Delta)$ and $G^+\Delta = \Delta G^+$.\\ii) let $g\in G$; there exists $g_1\in G^+$ and $n\in \mathbb{N}$ such that $g = g_1\Delta^{-n}$.\\ iii) let $g\in G$; there exists a unique couple $(a,b)$ of elements of $G^+$ such that $g = a^{-1}b$ and $a\land_L b = 1$. Furthermore if $cg$ is in $G^+$ for some $c$ in $G^+$ then $c\succ a$.
\end{Lem}   
The couple $(a,b)$ is called the left normal form of $g$ in $G$. In the same way, one can define the right normal form of $g$.\\

From the axioms of a Garside monoid, we can construct a left normal form on $G^+$ which generalises the classical (left) normal form on the Artin-Tits monoids~: There exists a well-defined function  $\alpha_L~:G^+\to D(\Delta)$ which associates to each element $w$ of $G^+$ the greatest element of $D(\Delta)$, for left-divisibility, which left-divides $w$; in particular we have $\alpha_L(w) = w\land_L\Delta$. For every element $w$ of $G^+-\{1\}$ we obtain a normal form $(w_1,\cdots, w_n)$ which is defined by $w = w_1\cdots w_n$ with $w_n\neq 1$ and $\alpha_L(w_i\cdots w_n) = w_i$ for every $i$. This normal form is ``local'', that is to say that  $(w_1,\cdots, w_n)$ is in normal form if and only if for every $i$, the couple $(w_i,w_{i+1})$ is in normal form. Furthermore,  one has $\alpha_L(wz) = \alpha_L(w\alpha(z))$.
We can define in the same way a right normal form and an associated function $\alpha_R$. Most of the time we will only use the left normal form~; so unless otherwise stated, the normal form will mean the left normal form; in particular, we write $\alpha$ for $\alpha_L$. The functions $\alpha_L$ and $\alpha_R$ can be defined in the same way and  without difficulties in the more general context of quasi-Garside group (see \cite{Dig}).

From now on, we fix a Garside system $(G,D(\Delta))$ associated to a   positive Garside  system $(G^+,D(\Delta))$. We denote by $S$ the set of atoms of $G^+$. We are concern with the Garside groups but we develop the theory in the context of quasi-Garside groups since most of the proof works in this wider context. Note that most of the following theory could be done in that context of pre-Garside groups (see Section \ref{sectionssgpatom}); except that the notion of parabolic subgroup is no more as natural as in context of the Garside groups because a pre-Garside group cannot have a Garside element) In that case the interesting notion is probably the notion of atomic Garside groups. We prefer to restrict ourself to Garside systems because in the context of pre-Garside groups, it is not clear how to go from the monoid to the group---pre-Garside monoids do not verify in general the Ore relations---and, in particular fundamental questions, such as to solve the word problem, remain open and should be answered before to try to go farther.    
\section{Classical Parabolic subgroups and two other relevant examples}
\subsection{The classical case of Artin-Tits groups} \label{sec.cas.class} We refer for Appendix A, for a definition of Artin-Tits groups and for the notation.  A subgroup $A_X$ ({resp.} submonoid $A_X^+$) of an Artin-Tits group $A$ ({resp.} monoid $A^+$) generated by a subset $X$ of $S$ is called a standard parabolic subgroup ({resp.} parabolic submonoid). Every subgroup of $A$ conjugate to a standard parabolic subgroup is called a parabolic subgroup. Van Der Lek showed in \cite{VdL} that $(A_X,X)$ is Artin-Tits system canonically isomorphic to the Artin-Tits system associated to the matrix $(m_{s,t})_{s,t\in X}$. The Garside element of $A_X$  is balanced in $A$ and is the lcm of $X$. It is not difficult to see that there is a 1-1 correspondence between the balanced minimal elements of $A^+$ and the standard parabolic subgroups of $A$.  Furthermore, we have :
\begin{Prop} Let $A^+$ be a Artin-Tits monoid; then, \\(a) every parabolic submonoid of $A^+$ is closed by gcd and lcm;\\(b) every element of a parabolic submonoid has the same normal form considered as an element of $A^+$ or considered as an element of the submonoid;\\(c) the intersection of two parabolic submonoids $A^+_X$ and $A^+_Y$  is the parabolic submonoid $A_{X\cap Y}^+$.\end{Prop}
\begin{Prop} Let $A$ be a Artin-Tits group; then,\\ (d) For every subset $X$ of $S$, we have $A_X\cap A^+ = A_X^+$ ;\\(e) if $A_X$ is a standard parabolic subgroup of $A$ and $x$ belongs to $A_X$, then  the fractional decomposition  of $x$ in $A_X$  coincides with its fractional decomposition in $A$.\\ (f) the intersection of two standard parabolic subgroups $A_X$ and $A_Y$ is the standard  parabolic subgroup $A_{X\cap Y}$.\end{Prop}
Furthermore, as recall in Appendix \ref{lab.appendLCM}, part of these results remains true for the subgroups which are the image of an LCM-homomorphism. These results explain why we think that Proposition 1,2 and 3 show that our definition of a parabolic subgroup, and  more generally of a Garside subgroup, is the good ones in the context of Garside groups. The second example below is probably another strong argument for the correctness of our definition.\\  
 
\subsection{Two examples} The aim of these examples is to help the reader to follow us in the more technical Section 3. Hence we are going to claim here without clear argument that some subgroups are (atomic) Garside subgroups or parabolic subgroups or, neither one nor the other.     
\subsubsection{First example}\label{labexemple1} Consider the set $S = \{x,y,z\}$ and the group $G$ with presentation $$G = \langle S \mid x^2  = y^2 ; xz = zx; yz = zy\rangle.$$ Then the group $G$ is the direct product of its two subgroups $G_{x,y}$ and $G_z$  generated by $\{x,y\}$ and $z$ respectively. $(G_{x,y},\{1,x,y,x^2\})$ and $(G_z,\{1,z\})$ are Garside systems (see \cite{DeP} for instance). It follows that $(G,D(x^2z))$ is a Garside system where $D(x^2z) = \{1;x;y;x^2;xz;yz;x^2z\}$. That Garside system provides us with various relevant examples on what should be called a Garside ({\it resp.} standard parabolic) subgroup of a Garside group and what should not be called a Garside ({\it resp.} standard parabolic)  subgroup.\\

The subgroups $G_{x,y}$ and $G_z$ are Garside subgroups with respective Garside elements $x^2$ and $z$; they are even two standard parabolic subgroups of $G$.\\ Consider $xz$ and the subgroup $G_{x,z}$  of $G$ generated by the support $\{x,z\}$ of $xz$. Then $(G_{x,z},\{1,x,z,xz\})$ is a Garside system. But $G_{x,z}$ is not acceptable as a standard parabolic subgroup of $G$ with $xz$ for Garside element because $x^2$ is in normal form in $G$ but its normal form in $G_{x,z}$ is $(x,x)$; we ask that a positive element of a standard parabolic (Garside) subgroup has the same normal form in the group and in the subgroup. This example shows that we cannot expect a 1-1 correspondence between the balanced minimal elements and the standard parabolic subgroups. But $G_{x,z}$ is an atomic Garside subgroup of $G$ with Garside element $x^2z$. \\Now consider $G_{y,z}$ the subgroup of $G$ generated by $\{y,z\}$; then $G_{y,z}$ is another atomic Garside subgroup of $G$ with Garside element $y^2z = x^2z$. The intersection of $G_{x,z}$ and $ G_{y,z}$ is the subgroup $G_{x^2,z}$ generated by $x^2$ and $z$. $G_{x^2,z}$ is not an atomic Garside subgroup because it is not generated by atoms of $G$. Nevertheless, $G_{x^2,z}$ is a Garside subgroup of $G$ with $x^2$ and $z$ for atoms and with $x^2z$ for Garside element. Hence we cannot expect the intersection of two atomic Garside subgroups to be in general an atomic Garside subgroup.
\subsubsection{Second example: The Birman-Ko-Lee presentation}
Consider the braid group $\mathcal{B}_{n}$ on $n$ strings. Recall that this group has two well-known presentations, namely the classical presentation and the dual presentation (see Appendix A for the definitions and the notation) and that each of them gives to $\mathcal{B}_{n}$ a structure of Garside group.\\ 

Consider the dual presentation of $\mathcal{B}_{n}$. Then the standard parabolic subgroups of the classical presentation remains standard parabolic subgroups for the dual presentation. For instance, $G_{s_1,s_2}$ is equal to $G_{a_{21},a_{32}}$ and its Garside element is $a_{21}a_{32}$ for the dual presentation. But we obtain other standard parabolic subgroups. For instance for $n = 4$ , $G_{a_{31},a_{21}}$ is a standard parabolic subgroup of $\mathcal{B}_{n}$ for the dual presentation. Remark that $G_{a_{13},a_{24}}$ is not a standard parabolic subgroup (neither a Garside subgroup) because the least common multiple of $a_{13}$ and $a_{24}$ in $\mathcal{B}_{n}$ for left-divisibility is not in $G_{a_{13},a_{24}}$. Therefore, the first natural idea in order to define a standard parabolic subgroup of a Garside group, that is to consider every subgroup generated by a subset of atoms of the Garside group, does not work. It is not difficult to verify that for $n= 4$, all the standard parabolic subgroups of the dual presentation are parabolic subgroups of the classical presentation (that is a classical standard parabolic subgroup or a subgroup conjugate to a standard parabolic subgroup). We conjecture that this fact is true for every Artin-Tits groups of spherical type. In that case it will make the notion of parabolic subgroup more natural than the notion of standard parabolic subgroup since the parabolic subgroups do not depend on the presentation. We note that this conjecture is no more true for the Affine braid group because all the standard parabolic subgroups of the classical presentation of that group are of spherical type whereas in the dual case, some standard parabolic subgroups are of Affine type (see \cite{Dig} where they are called quasi-parabolic subgroups). This implies that they cannot be classical parabolic subgroups (consider the center for instance).    
\section{Garside and standard parabolic subgroups}
\subsection {The Garside subgroups}
\begin{Def}[Garside subgroups]\ \label{defgarssbgp}\\
i) Let $(G^+,D(\Delta))$ be a   positive Garside  system. Let $H^+$ be a submonoid of $G^+$. We say that $H^+$ is a Garside submonoid of $G^+$ if :\\(1)
 $H^+$ is a sublattice of $G^+$ for left-divisibility and for right-divisibility.\\(2) $\alpha_R(H^+) \subset H^+$ and $\alpha_L(H^+) \subset H^+$.\\ii) Let $(G,D(\Delta))$ be a Garside system and $H$ be a subgroup of $G$. We say that $H$ is a Garside subgroup of $G$ if $H^+ = H\cap G^+$ is a Garside submonoid of $G^+$ and generates $H$.
\end{Def}
Note that by definition, a Garside submonoid is closed by l.c.m. and by g.c.d. As a consequence we have 
\begin{Lem} Let $(G^+,D(\Delta))$ be a   positive Garside  system and $H^+$ be a Garside submonoid of $G^+$. Let $h_1,h_2$ be in $H^+$; then $h_1$ left-divides $h_2$ in $G^+$ if and only if $h_1$ left-divides $h_2$ in $G^+$. The same result is true if we replace  left-divisibility by right-divisibility.
\end{Lem}
\begin{proof}
\end{proof}
With the above definition we have :
\begin{The}\label{labTh1}
Let $(G,D(\Delta))$ be a Garside system and $H$ be a Garside subgroup of $G$. Set $H^+ = H \cap G^+$;
then,\\(i) $H^+\cap D(\Delta)$ is a sublattice of $G^+$ for left-divisibility and for right-divisibility. Furthermore, the left upper-bound and the right upper-bound are equal.\\(ii) $(H^+, D(\Delta)\cap H^+)$ is a positive Garside  system; its Garside element is the upper-bound of $H^+\cap D(\Delta)$. Furthermore, if $h\in H^+$ and $(h_1,\cdots,h_n)$ is its normal form in $G^+$, then $h_i\in H^+$ for every $i$. In particular $(h_1,\cdots,h_n)$ is the normal form of $h$ in $H^+$.\\
(iii) $H$ is isomorphic to the group of fractions of $H^+$. Furthermore, if $h = a^{-1}b$ is the left normal form of $h$ in $G$, then $a$ and $b$ are in $H^+$. In particular, $h = a^{-1}b$, is also the left normal form of $h$ in $H$. The same is true for the right normal form.   
\end{The}
\begin{proof} (i) The intersection of two sublattices is a sublattice. Then $H^+\cap D(\Delta)$ is a sublattice of $G^+$ for left-divisibility and for right-divisibility. Denote by $\delta_L$ and $\delta_R$ the left upper-bound and the right upper bound of $H^+\cap D(\Delta)$ respectively. Then $\delta_R$ left-divides $\delta_L$  and $\delta_L$ right-divides $\delta_R$. By the remark on the Noetherian hypothesis which follows Definition \ref{defgarsgp}, $\delta_R = \delta_L = \delta$.\\(ii) The element $\delta$ is balanced in $H^+$ because $D_{H^+}(\delta) = D(\delta)\cap H^+$ is both the set of its left-divisors and of its right-divisors in $H^+$. In particular the axiom (b') of the definition of positive Garside systems is true. Now consider $g$ in $H^+$ and $(h_1,\cdots,h_n)$ its normal form in  $G^+$. The element $h_1 = \alpha(g)$ is in $H^+\cap D(\Delta) = D_{H^+}(\delta)$. Since $H^+$ is a sublattice of $G^+$, we can write  $g = h_1z$ with $z$ in $H^+$. By cancellativity, we get that $z = h_2\cdots h_n$. Hence $h_2\cdots h_n$ is in $H^+$ with $(h_2,\cdots,h_n)$ for normal form. By induction on $n$, we get that $h_2,\cdots, h_n$ are in $D_{H^+}(\delta)$. As a consequence, $D_{H^+}(\delta)$ is a generating set for $H^+$ and $(H^+,D(\delta))$ is a positive Garside system. Furthermore the left normal form in $H^+$ of the above element $g$ is $(h_1,\cdots, h_n)$: consider $(h'_A,\cdots, h'_m)$ its left normal form in $H^+$. On the one hand $h'_1$ is in $D(\Delta)$ and then $h'_1$ divides $h_1$; and the other hand, $h'_1$ is in $D_{H^+}(\delta)$ and then $h_1$ divides $h'_1$. It follows that $h_1 = h'_1$ and we conclude by induction on $n$.\\(iii) Consider $g$ in $H$ and $a^{-1}b$ its normal form in $H$.  We have $a\land_L b = 1$ in $H^+$ and because $H^+$ is a sublattice, we have also $a\land_L b = 1$ in $G^+$. Thus by unicity of the left normal form, $a^{-1}b$ is the normal form of $g$ in $G$. The isomorphism follows.  \end{proof}

\begin{Prop} \label{labprop1}Let $(G,D(\Delta))$ be a Garside system and $H$, $K$ be two Garside subgroups of $G$. Then $H\cap K$ is a Garside subgroup of $G$ and $(H\cap K)\cap G^+ = H^+\cap K^+$.  
\end{Prop}
\begin{proof} It is clear that $(H\cap K)\cap G^+$ is equal to $H^+\cap K^+$ and verifies the two properties which define a Garside submonoid. Now let $g$ be in $H\cap K$ and consider $a^{-1}b$ its normal form in $G$. Then $a$ and $b$ are in $H^+$ and in $K^+$. Hence,  $H^+\cap K^+$ generates $H\cap K$. 
\end{proof}

It is obvious that the Garside element of the Garside subgroup $H\cap K$ divides the left (see Lemma \ref{lemsuivant}) g.c.d of the Garside elements of $H$ and $K$. We do not see why this two elements should be equal but we cannot find a counterexample.\\

Note that Theorem \ref{labTh1} and Proposition \ref{labprop1} prove Proposition \ref{labpropintro1} and \ref{labpropintro2}.\\
\subsubsection{Recognising a Garside subgroup}\label{secrecogn} 
Now, the definition of a Garside subgroup is not very practical. In particular, it is natural to address the following question : Given $(G,D(\Delta))$ a Garside system and $X$ a sublattice of $D(\Delta)$---that is a subset which is a sublattice for both left and right divisibility(consequently the two upper-bounds are equal)---which contains $1$. How can I know if $G_X$, the subgroup of $G$ generated by $X$, is a Garside subgroup of $G$ with $X$ as set of minimals ?  The following proposition answers that question. 
\begin{Prop} Let $(G,D(\Delta))$ be a Garside system and $X$ a sublattice of $D(\Delta)$ which contains $1$. Set $X^2 = \{xy\mid x,y\in X\}$; then $G_X$ is a Garside subgroup with $X$ for set of minimals if and only if $\alpha_L(X^2) = X$ and $\alpha_R(X^2) = X$.
\end{Prop} 
\begin{proof} Set $G_X^+ = G^+\cap G_X$. If $G_X$ is a Garside subgroup with $X$ for set of minimals then $\alpha(X^2) = X$ for left-divisibility and for right-divisibility.. Conversely, assume that $\alpha(X^2) = X$ for left-divisibility and for right-divisibility. It is clear that $(G_X^+, X)$ is a positive Garside system. Let us first shows that the left normal form $(h_1,\cdots,h_n)$ in $G_X^+$ of $h \in G_X^+$ is also its normal form in $G^+$; that is $\alpha(h_ih_{i+1}) = h_i$ for every $i$. We prove it by induction on $n$. For $n = 1$ it is obvious. Since $(h_2,\cdots, h_n)$ is in left normal form in $G^+_X$, we have by induction hypothesis that $(h_2,\cdots, h_n)$ is in left normal form in $G^+$. So it remains to prove that $\alpha(h_1h_2) = h_1$. On the one hand, $h_1$ is in $D(\Delta)$, then $h_1$ left-divides $\alpha(h_1h_2)$ in $G^+$. On the other hand, $\alpha(h_1h_2)$ is in $X$ which is a sublattice of $D(\Delta)$.  Then $h_1$ left-divides  $\alpha(h_1h_2)$ in $G_X^+$. Write $\alpha(h_1h_2) = h_1z$ with $z$ in $X$. But  $h_1z$ left-divides $h_1h_2$ in $G^+$ and $z$ left-divides $h_2$ in $G^+$. As before we get that $z$ left-divides $h_2$ in $G^+_X$. Thus $h_1z$ left-divides $h_1h_2$ in $G^+_X$. since $(h_1,h_2)$ is in left normal form in $G^+_X$, we get $z = 1$ and $h_1 = \alpha(h_1h_2)$ in $G^+$. Hence an element of $G_X^+$ has the same left normal form in $G^+$ and in $G^+_X$. Of course the same is true for the right normal form. Now we prove that $G^+_X$ is a sublattice of $G^+$ as in the proof of Theorem 2.10 of \cite{god_agt} (the fact that $X$ is a sublattice replace here Lemma 2.9 of \cite{god_agt} in the that proof).   \end{proof} 

Now if we take every no-empty subset $X$ of $D(\Delta)$ in order to know if the subgroup $G_X$ of $G$ generated  by $X$ is a Garside subgroup of $G$, we have  to consider the set of minimals of $G^+$ which can be written as a product of elements of $X$; to verify  firstly that this set is a sublattice of $D(\Delta)$   
and secondly the above criterion. When $G$ is a Garside group, that is when $D(\Delta)$ is finite, this method gives us a finite set of properties to verify in order to know if $G_X$ is a Garside subgroup of $G$.

\subsection{The parabolic subgroups}
Let $(G,D(\Delta))$ be a Garside system and $\delta$ be a balanced element of $D(\Delta)$. We denote by $G_\delta$ the subgroup of $G$ generated by $supp(\delta)$, and we set $G^+_\delta = G_\delta\cap G^+$.
\begin{Prop} Let $(G,D(\Delta))$ be a Garside system. Let $\delta$ be a balanced element of $D(\Delta)$; then $(G_\delta,D(\delta))$ is a Garside system with $supp(\delta)$ for set of atoms. Furthermore, $G_\delta$ is a Garside subgroup with $\delta$ for Garside element if and only if $D(\delta) = D(\Delta)\cap G^+_\delta$. 
\end{Prop}
\begin{proof} It is clear that assertions (a), (b') and (c) of Definition \ref{defgarsgp} are true; the first claim of the proposition follows. Now, by the general results on Garside subgroups, if $G_\delta$ is a Garside subgroup with $\delta$ for Garside element then $D(\delta) = D(\Delta)\cap G^+_\delta$.
Conversely, assume that $D(\delta) = D(\delta)\cap G^+_\Delta$. So $D(\delta)$ 
is a sublattice for left-divisibility and for right-divisibility. Let $x$, $y$ be in 
$D(\delta)$ such that $(x,y)$ is in normal form in $G^+_\delta$ and assume that $\alpha(xy)$ in $G^+$ is not equal to $x$. Then there exists an atom $s$ of $G^+$ such that $xs$ is in $D(\Delta)$ and left-divides $xy$. But in that case $s$ left-divides $y$ (in $G^+$ and in $G^+_\delta$) and as a consequence is in $supp(\delta)$. Hence $xs$ is in $G^+_\delta\cap D(\Delta) = D(\delta)$, and a contradiction. Then $\alpha(xy)$ is in  $D(\delta)$ for every $x$, $y$ in $D(\delta)$. The have the same result for the right normal form. Applying the result of Subsection \ref{secrecogn}, we are done. \end{proof}
\begin{Def}[Parabolic subgroups]\ \\ Let $(G,D(\Delta))$ be a Garside system.\\1) let $\delta$ be a balanced element of $D(\Delta)$. Denote by $G_\delta$ the subgroup generated by $supp(\delta)$. We say that $G_\delta$ is a standard parabolic subgroup if $D(\delta) = D(\Delta)\cap G^+_\delta$. In that case, $G^+_\delta$ is called a parabolic submonoid of $G^+$ \\2) let $(G,D(\Delta))$ be a Garside system and $H$ be a subgroup of $G$. We say that $H$ is a parabolic subgroup of $G$ if $H$ is conjugate to a standard parabolic subgroup. 
\end{Def}

Since we only deal with standard parabolic subgroups, we will say in the following {\it parabolic subgroup} for {\it standard parabolic subgroup}. \begin{Lem}\label{axiomenplus} Let $(G^+,D(\Delta))$ be a positive Garside system and $G^+_\delta$ be a parabolic submonoid of $G^+$. Then $G^+_\delta$ is closed by left-divisibility and by right-divisibility. \end{Lem}
\begin{proof} By symmetry, it is enough to prove the result for left-divisibility. So consider $g,h$ in $G^+$ with $h$ in $G^+_\delta$ and assume that $g$ left-divides $h$. Denote by $(h_1,\cdots, h_n)$ and $(g_1,\cdots, g_m)$ the left normal forms of $h$ and $g$ respectively. We have that $g_1$ left-divides $h_1$ because $g$ left-divides $h$. Then $g_1$ is in $D(\delta)$ and in $G^+_\delta $. Furthermore, if we write $h_1 = g_1h'_1$, with $h'_1$ in $D(\delta)$, then $g_2\cdots g_m$ left-divides  $h'_1h_2\cdots h_n$. Since $(g_2,\cdots, g_m)$ is in normal form, we  can apply an induction argument on $m$ to conclude. 
\end{proof}
\begin{Lem} \label{lemsuivant}Let $(G^+,D(\Delta))$ be a positive Garside system and $\delta$, $\tau$ be two balanced elements of $G^+$; then $\delta\land_L\tau$ and $\delta\land_R\tau$ are balanced and equal in $G^+$. 
\end{Lem}
\begin{proof} By definition $\delta\land_L\tau$ left-divides both $\delta$ and $\tau$. Since these two elements are balanced, we get that $\delta\land_L\tau$ right-divides both $\delta$ and $\tau$ too. Then it right-divides $\delta\land_R\tau$. But by the same arguments we have also that  $\delta\land_R\tau$ left-divides $\delta\land_L\tau$. By the Noetherian hypothesis, it implies that  $\delta\land_L\tau$ and $\delta\land_R\tau$ equal and consequently balanced. 
\end{proof}
Since $\delta\land_L\tau = \delta\land_R\tau$ when $\delta$ and $\tau$ are balanced, we will write $\delta\land\tau $ for that element.  
\begin{Prop} Let $(G,D(\Delta))$ be a Garside system and $G_\delta$, $G_\tau$ be two parabolic subgroups of $G$. Then $G_\delta\cap G_\tau = G_{\delta\land\tau}$ and is a parabolic subgroup of $G$.  
\end{Prop}
\begin{proof} By Proposition \ref{labprop1}, $G_\delta\cap G_\tau$ is a Garside subgroup with $G^+_\delta\cap G^+_\tau\cap D(\Delta) = D(\delta)\cap D(\tau) = D(\delta\land \tau)$ for set of minimals. 
\end{proof}
\subsection{The atomic Garside subgroups}\label{labgeneralisation}\label{sectionssgpatom}
\subsection{Definition}
\begin{Def} 1) Let $(G^+,D(\Delta))$ be a positive Garside system and $H^+$ be a Garside submonoid of $G^+$. Denote by $S$ and $S_{H^+}$ the sets of atoms of $G^+$ and $H^+$ respectively. We say that $H^+$ is an atomic Garside submonoid of $G^+$ when $S_{H^+}$ is a subset of $S$. \\2) Let $(G,D(\Delta))$ be a Garside system and $H$ be a Garside subgroup of $G$. We say that $H$ is an atomic Garside subgroup of $G$ if $H^+ = H\cap G^+$ is an atomic Garside submonoid of $G^+$.\end{Def}
Note that by definition, a (standard) parabolic subgroup is an atomic Garside subgroup.\\

As we have seen in Example \ref{labexemple1}, the intersection of two atomic Garside subgroups is not in general an atomic Garside subgroup. So we address the following question: what can we say of the intersection of two atomic Garside subgroups ?   
Let us finish this section by the easy following lemma.
\begin{Lem} Let $(G^+,D(\Delta))$ be a positive Garside system. Let $\delta$ be a balanced element of $G^+$. Then,\\(i) $D(\delta)\delta = \delta D(\delta)$ and $supp(\delta)\delta = \delta supp(\delta)$. \\(ii) if $H^+$ is an atomic Garside submonoid of $G^+$ with $\delta$ for Garside element then the set of atoms of $H^+$ is an union of conjugacy classes of $supp(\delta)$ for the conjugacy by $\delta$. \end{Lem} 
\begin{proof} Let $x$ be in $D(\delta)$. Then there exists $y$ in $D(\delta)$ such that $xy = \delta$. Since $\delta$ is balanced, there exists $z$ in $D(\delta)$ such that $yz = \delta$. We get $x\delta = xyz = \delta z$. Thus $D(\delta)\delta \subset \delta D(\delta)$. By symmetry, we have also the other inclusion and then the equality. Now it is easy to deduce from $D(\delta)\delta = \delta D(\delta)$  that $supp(\delta)\delta = \delta supp(\delta)$. Finally, if $x$ is an atom of $H^+$ and $x\delta = \delta z$, we get that $z$ is in $H^+$ because  $H^+$ is a sublattice and that both $x\delta$ and $\delta$ are in $H^+$. \end{proof}
\subsubsection{Pre-Garside monoid} 
An alternative but equivalent approach to Garside group was introduced in \cite{BDM}. We refer to that paper for the notion of pre-monoid and pre-Garside monoid. A pre-Garside group is the group of fractions of a pre-Garside monoid. In \cite{BDM} the authors gave an easy characterisation of pre-Garside monoids which are Garside monoids. Conversely, every Garside monoid can be seen as a pre-Garside monoid. In the context of pre-Garside groups, we do not have in general a Garside element. The main examples of pre-Garside monoid are the Artin-Tits monoid, not necessarily of spherical type. In that context, the notion of atomic Garside sub-monoid is probably more natural than our notion of standard parabolic subgroups since we do not refer to a Garside element. Furthermore one can see that in the case of Artin-Tits monoid, it corresponds precisely with the classical notion of parabolic submonoid. Consequently we address the following question: What is the good definition of a standard parabolic submonoid ({\it resp.} subgroup) in the context of pre-Garside monoids({\it resp.} groups)? We mean, which axiom(s) should we add to the axioms of the definition of an atomic Garside subgroup if we want that the family of (classical) standard parabolic subgroups of every Artin-Tits groups becomes the family of (the new) standard parabolic subgroups associated to the classical presentation, and that, in the general context of Garside group, the family of standard parabolic subgroups is closed by intersection?  Perhaps, we should consider the property proved  in Lemma \ref{axiomenplus} as the extra axiom. Since it is not the topic of that paper, we do not go farther on that question.\\

\noindent{\bf Acknowledgements:} I thank Patrick Dehornoy and Fran\c cois Digne for useful remarks and comments which help me to improve the general exposition of my results.   

\appendix  
\section{Artin-Tits groups.}
In this Appendix, we recall the definition of the Artin-Tits groups and their relevant properties regarding the present paper. \subsection{Definition.}
Let $S$ be a finite set and $M= (m_{s,t})_{s,t\in S}$ a symmetric matrix  with $m_{s,s} = 1$
for $s\in S$ and $m_{s,t} \in \mathbb{N}-\{0;1\}\cup \{\infty\}$
for $s\not=t$ in $S$.  The Artin-Tits system
associated to $M$ is the pair $(A,S)$   where $A$ is the group defined by the presentation
$$A = \langle S|\underbrace{sts\cdots}_{m_{s,t}\ terms} =
\underbrace{tst\cdots}_{m_{s,t}\ terms}~;\ \forall s,t\in S, s\not= t\
and\ m_{s,t}\not=\infty \rangle.$$ The group $A$ is called an Artin-Tits group
and relations $\underbrace{sts\cdots}_{m_{s,t}\ terms} =
\underbrace{tst\cdots}_{m_{s,t}\ terms}$ are called braid relations. For
instance, if
$S = \{s_1,\cdots, s_n\}$ with $m_{s_i,s_j}= 3$ for $|i-j| =1$ and
$m_{s_i,s_j} = 2$ otherwise, then the associated Artin-Tits group is the braid group $\mathcal{B}_n$. We
denote by $A^+$ the submonoid of $A$ generated by  $S$. This monoid
$A^+$ has the same presentation as the group $A$, considered as a monoid
 presentation (\cite{Par4}); its elements are called the positive elements of $A$. When we add the relations $s^2 = 1$ for $s\in S$ to the presentation of $A$ we obtain the Coxeter group $W$ associated to $A$. We say that $A$ (or simply S) is of spherical type if $W$ is finite. In that case, $S$ has a left-least common multiple $\Delta_S$ in $A^+$ which is also its right-least common multiple. Denote by $i:A\to W $ the canonical surjection. Then $i$ has a canonical section (see \cite{Bou} chapter 4)~; furthermore an element of the image $D(\Delta_S)$ of this section is characterised by the fact every word which represents it is square free. That is to say it does not contain a square of an element of $S$ in its writing. Furthermore when $A$ is of spherical type, then $(A,D(\Delta_S))$ is a Garside system and $i(\Delta_S)$ is the longest element of $W$.
\subsection{LCM-homomorphism.}\label{lab.appendLCM}
\begin{Def}[\cite{Cri},\cite{god_agt} Definition 0.1]\label{defLcm}
Let $(A,S)$ and $(B,T)$ two Artin-Tits systems and $p$ an
 application from $S$ in $ \mathcal{P}(T)-\{\emptyset\}$, the set of non-empty subsets of  $T$, such that \\ {\bf(L0)} if $s\not=t\in S$ then $p(s)$ and $p(t)$ are disjoint;\\ {\bf(L1)} for every $s\in S$, $p(s)$ is of spherical type;\\
{\bf (L2)}  if $s\not=t\in S$ with $m_{s,t}\not= \infty$, one has\\
 \begin{center} $\underbrace{\Delta_{p(s)}\Delta_{p(t)}\cdots}_{m_{s,t}\ terms}  = \underbrace{\Delta_{p(t)}\Delta_{p(s)}\cdots}_{m_{s,t}\ terms}  = \Delta_{p(s)}\lor_L \Delta_{p(t)}$ in $B^+$\end{center} {\bf(L3)} if $s\not=t\in S$ with $m_{s,t}= \infty$, then \begin{center}$\left\{\begin {array}{l} \forall u\in p(s),\ \{u\}\cup p(t)\textrm{ is not of spherical type},\\ \forall u\in p(t),\ \{u\}\cup p(s)\textrm{ is not of spherical type.}\end{array}\right.$ \end{center} One can then define an homomorphism $\varphi_p$ from $A$ to $B$ by setting $\varphi_p(s)= \Delta_{p(s)}$ for $s\in S$. Such a morphism is called an LCM-homomorphism.
\end{Def}
If we focus on Artin-Tits groups of spherical type then the axiom (L3) can be cancelled since the case $m_{s,t} = \infty$ never happens. 
\begin{The}[\cite{Cri} Theorem 1.3 and 1.4]\label{Thdeb2} Let $(A,S)$ and $(B,T)$ be two Artin-Tits systems of spherical type and $\varphi_p~:~A\to B$ an LCM-homomorphism~; then $\varphi_p$ is injective. \end{The}
\begin{The}\label{cor2deb} Let $(A,S)$ and $(B,T)$ be two Artin-Tits systems and $\varphi_p: A\to B$ an
 LCM-homomorphism~; then\\(i) {\bf (\cite{god_agt} Theorem 2.10)} 
 $\varphi_p$ respects the normal forms~: if $(g_1,\cdots,g_n)$ is the normal form of $g\in A^+$ then $(\varphi_p(g_1),\cdots,\varphi_p(g_n))$ is the normal form of $\varphi_p(g)$ in $B^+$.\\(ii) {\bf (\cite{Cri2} Theorem  8, \cite{god_agt} Theorem 2.10)} $\varphi_p$ respects the lcm and the gcd~: if $g,h\in A^+$ then $\varphi_p(g)\lor_L \varphi_p(h) ) = \varphi_p(g\lor_L h)$ and $\varphi_p(g)\land_L \varphi_p(h) ) = \varphi_p(g\land_L h)$ ~; the same is true if we replace $\prec$ by $\succ$. \end{The} 
\subsection{The Birman-Ko-Lee presentation}
In \cite{BKL}, Birman, Ko and Lee proposed an alternative presentation of the braid group on n strings. This point of view was generalised by Bessis in \cite{Bes} to every Artin-Tits groups of spherical type. 
\begin{Prop} The braid group $\mathcal{B}_n$ has the following presentation :
$$\displaylines{\mathcal{B}_n  = \left\langle\  a_{ts}\  ; n\geq t\geq s\geq 1\mid\begin{array}{ll} \ \\ \ \end{array}\right. \hfill\cr\hfill\left.\begin{array}{ll} a_{ts}a_{rq} = a_{rq}a_{ts}&\textrm{when } (t-r)(t-q)(s-r)(s-q) > 0 \\ a_{ts}a_{sr} = a_{sr}a_{tr} = a_{tr}a_{ts}&\textrm{when } t>s>r\end{array}\right\rangle}$$\end{Prop} 

That presentation is called the dual presentation of the braid group whereas the one given in A.1 is called the classical presentation. The generator $s_i$ of the classical presentation is equal to the generator $a_{(i+1)i}$ of the dual presentation. \\
Furthermore if $\mathcal{B}^{BKL+}_n$ is the submonoid of $\mathcal{B}_n$ generated by the $a_{ts}$. Then $\delta = a_{n(n-1)}\cdots a_{21}$ is balanced in $\mathcal{B}^{BKL+}_n$ and $(\mathcal{B}_n, D(\delta))$ is a Garside system (see \cite{DeP}).

\end{document}